\nonstopmode
\documentclass[a4paper, 10pt, reqno]{amsart}
\usepackage{latexsym}
\usepackage{fancyhdr}
\usepackage{amsmath, amssymb}
\usepackage[ansinew]{inputenc}
\usepackage[all]{xy}
\usepackage{verbatim}
\usepackage{graphicx}
\theoremstyle{plain}
\newtheorem*{lemma*}{Lemma}
\newtheorem{lemma}[subsection]{Lemma}
\newtheorem*{theorem*}{Theorem}
\newtheorem{theorem}[subsection]{Theorem}
\newtheorem*{proposition*}{Proposition}
\newtheorem{proposition}[subsection]{Proposition}
\newtheorem*{corollary*}{Corollary}

\newtheorem*{result*}{Result}
\newtheorem{result}[subsection]{Result}
\theoremstyle{definition}
\newtheorem*{definition*}{Definition}

\newtheorem*{example*}{Example}
\newtheorem{example}[subsection]{Example}
\theoremstyle{remark}
\newtheorem*{remark*}{Remark}
\newtheorem{remark}[subsection]{Remark}
\newtheorem*{remarks*}{Remarks}
\newtheorem{remarks}[subsection]{Remarks}

\numberwithin{equation}{section}

\def\al{\alpha}

\def\si{\sigma}

\def\ta{\tau}

\def\om{\omega}

\def\De{\Delta}

\def\Om{\Omega}

\def\R{\mathbb{R}}

\def\aut{\mathrm{Aut}}

\def\id{\mathrm{id}}

\def\reg{\mathrm{reg}}

\def\<{\langle}
\def\>{\rangle}
\def\Hyp{\mathrm{Hyp}}
\def\sing{\mathrm{sing}}

\def\orb{\mathrm{orb}}
\def\amb{\mathrm{amb}}

\let\on=\operatorname
\def\sr#1%
{\ifmmode{}^\dagger\else${}^\dagger$\fi\ifvmode
\vbox to 0pt{\vss
 \hbox to 0pt{\hskip\hsize\hskip1em
 \vbox{\hsize3cm\eightpoint\raggedright\pretolerance10000
 \noindent #1\hfill}\hss}\vss}\else
 \vadjust{\vbox to0pt{\vss%
 \hbox to 0pt{\hskip\hsize\hskip1em%
 \vbox{\hsize3cm\eightpoint\raggedright\pretolerance10000%
 \noindent #1\hfill}\hss}\vss}}\fi%
}

\title[Choosing roots of polynomials with symmetries smoothly]
{Choosing roots of polynomials with symmetries smoothly}

\author[M. Losik, A. Rainer]
{Mark Losik and Armin Rainer}

\address{M. 
Losik: Saratov State University, 
Astrakhanskaya, 83, 410026 Saratov, Russia}

\email{losikMV@info.sgu.ru}

\address{A. 
Rainer: Institut f\"ur Mathematik, Universit\"at Wien, 
Nordbergstrasse~15, A-1090 Wien, Austria}

\email{armin.rainer@univie.ac.at}

\begin{document}

\begin{abstract} 
The roots of a smooth curve of hyperbolic polynomials may not in 
general be parameterized smoothly, even not $C^{1,\al}$ for any $\al > 0$. 
A sufficient condition for the existence of 
a smooth parameterization is that no two of the increasingly ordered 
continuous roots meet of infinite order. We give refined sufficient 
conditions for smooth solvability if the polynomials have certain symmetries.
In general a $C^{3n}$ curve of hyperbolic polynomials of degree $n$ 
admits twice differentiable parameterizations of its roots. 
If the polynomials have certain symmetries we are able to weaken the 
assumptions in that statement. 
\end{abstract}

\thanks{The authors were supported by
`Fonds zur F\"orderung der wissenschaftlichen Forschung, Projekt P 17108-N04'}
\date{March 28, 2006}

\maketitle

\section{Introduction}
 
Consider a smooth curve of monic hyperbolic (i.e.\!\! all roots real) 
polynomials with fixed degree $n$:
\[
P(t)(x) = x^n - a_1(t) x^{n-1} +a_2(t) x^{n-2} - \cdots + (-1)^n a_n(t)
\qquad (t \in \mathbb R).
\]
Is it possible to find $n$ smooth functions $x_1(t),\ldots,x_n(t)$ 
which parameterize the roots of $P(t)$ for each $t$? 
It has been shown in \cite{rellich} that real analytic curves $P(t)$ 
allow real analytic parameterizations of its roots, and in \cite{roots}
that the roots of smooth curves $P(t)$ may be chosen smoothly if no 
two of the increasingly ordered continuous roots meet of infinite order.
In general, as shown in \cite{rootsII}, the roots of a $C^{3n}$ 
curve $P(t)$ of hyperbolic polynomials can be parameterized twice 
differentiable. That regularity of the roots is best possible: 
In general no $C^{1,\al}$ parameterizations of the roots for any $\al > 0$
exist which is shown by examples in \cite{roots}, \cite{BBCP}, and \cite{Glaeser}. 
Further references related to that topic are 
\cite{bronshtein}, \cite{Mandai}, and \cite{wakabayashi}.

The space $\Hyp^n$ 
of monic hyperbolic polynomials $P$ of fixed degree $n$ may 
be identified with a semialgebraic subset in $\R^n$, 
the coefficients of $P$ being the coordinates. 
Then $P(t)$ is a smooth curve in $\Hyp^n \subseteq \R^n$. 
If the curve $P(t)$ lies in some semialgebraic subset of 
$\Hyp^n$, then it is evident that in general the conditions which 
guarantee smooth parameterizations of the roots of $P(t)$ are weaker 
than those mentioned in the previous paragraph.
In the present paper we are going to study that phenomenon.

In section \ref{orbemb} we present a class of semialgebraic subsets 
in spaces of hyperbolic polynomials for which we are able to apply 
the described strategy. The construction of that class is based 
on results due to \cite{smithstrong}.

Our main goal is to investigate the problem of finding smooth roots 
of $P$ under the assumption that the polynomials $P(t)$ 
satisfy certain symmetries.
More precisely, we shall assume that the roots $x_1(t),\ldots,x_n(t)$ 
of $P(t)$ fulfill some linear relations, i.e., 
there is a linear subspace $U$ of $\mathbb R^n$ such that 
$(x_1(t),\ldots,x_n(t)) \in U$ for all $t$. 
Then the curve $P(t)$ lies in the semialgebraic subset $E(U)$ 
of the space of hyperbolic polynomials $\Hyp^n = E(\mathbb R^n) = 
\mathbb R^n/\on{S}_n$ of degree $n$, where $E = (E_1,\ldots,E_n)$ and 
$E_i$ denotes the $i$-th elementary symmetric function. 
The symmetries of the roots of $P(t)$ are represented by the action 
of the group $W$ on $U$ which is inherited from the action 
of the symmetric group $\on{S}_n$ on $\mathbb R^n$ by permuting the 
coordinates: 
\[
W=W(U):=N(U)/Z(U),
\]
where $N(U) := \{\tau \in \on{S}_n : \tau.U=U\}$ and 
$Z(U) := \{\tau \in \on{S}_n : \tau.x=x ~\mbox{for all}~ x \in U\}$.

Under the additional assumption that the restrictions $E_i|_U$, 
$1 \le i \le n$, generate the algebra $\mathbb R[U]^W$
of $W$-invariant polynomials on $U$, we will 
show that the conditions imposed on $P(t)$ in order to guarantee 
the existence of a smooth parameterization of its roots may be 
weakened.
These conditions will be formulated in terms of the two natural 
stratifications carried by $U$ and $E(U) = U/W$: 
the orbit type stratification with respect to $W$ and 
the restriction of the orbit type stratification 
with respect to $\on{S}_n$. 
The latter will be called ambient stratification.
See section \ref{orbamb}.
It will turn out (section \ref{polysym}) that we may find global smooth parameterizations 
of the roots of $P(t)$, provided that $P(t)$ is 
normally nonflat with respect to the orbit type stratification of 
$E(U) = U/W$ at any $t$. 
This condition is in general weaker than the condition 
found in \cite{roots}, since we prove in section \ref{orbamb} 
that normal nonflatness with respect to the ambient stratification 
implies normal nonflatness with respect to the orbit type stratification.
For a definition of `normally nonflat' see \ref{smoothlifts}. 

These improvements are essentially applications of the lifting 
problem tackled in \cite{rep-lift}. See also \cite{rep-lift3} and 
\cite{rep-lift4}. This generalization of the above problem studies 
the question whether it is possible to lift smoothly a smooth curve 
in the orbit space $V/G$ of an orthogonal finite dimensional 
representation of a compact Lie group $G$ into the representation space $V$. 
Here the orbit space $V/G$ is identified with the semialgebraic subset 
$\sigma(V)$ in $\mathbb R^n$ given by the image 
of the orbit map $\sigma = (\sigma_1,\ldots,\sigma_n): V \to \mathbb R^n$,
where $\sigma_1,\ldots,\sigma_n$ constitute a system of homogeneous 
generators of the algebra $\mathbb R[V]^G$ of $G$-invariant polynomials on 
$V$. See section \ref{pre} for details.

As mentioned before a $C^{3n}$ curve $P(t)$ of hyperbolic polynomials 
of degree $n$ allows twice differentiable parameterizations of its roots. 
Using results found for the general lifting problem in \cite{rep-lift4}, 
we are able to 
lower the degree of regularity in the assumption of that statement, 
if the polynomials $P(t)$ satisfy certain symmetries. 
See section \ref{diff}.  

A class of examples for which the described refinements 
apply will be constructed in section \ref{constr}.
For illustration we consider the case when $W$ is a finite reflection 
group in section \ref{refl}. 
Moreover, explicit examples will be treated.

The problem of finding regular roots of families of hyperbolic polynomials 
has relevance in the perturbation theory of selfadjoint operators 
(e.g. \cite{Kato}, \cite{perturb}, \cite{rellich}) and 
in the theory of partial differential equations for the well-posedness 
of hyperbolic Cauchy problems 
(e.g. \cite{BronshteinC}, \cite{hoermander2}).

\section{Preliminaries} \label{pre}

\subsection{Representations of compact Lie groups} \label{setting}

Let $G$ be a compact Lie group and let $\rho : G \rightarrow \on O(V)$ be an
orthogonal representation in a real finite dimensional Euclidean
vector space $V$ with inner product $\langle \quad \mid \quad \rangle$.
By a classical theorem of Hilbert and Nagata,
the algebra $\mathbb{R}[V]^{G}$ of invariant polynomials on $V$
is finitely generated.
So let $\sigma_1,\ldots,\sigma_n$ be a system of homogeneous generators
of $\mathbb{R}[V]^{G}$ of positive degrees $d_1,\ldots,d_n$.
Consider the \emph{orbit map}
$\sigma = (\sigma_1,\ldots,\sigma_n) : V \rightarrow \mathbb{R}^n$.
The image $\sigma(V)$ is a semialgebraic set in
$Z:=\{y \in \mathbb{R}^n : P(y) = 0 ~\mbox{for all}~ P \in I\}$
where $I$ is the ideal of relations between $\sigma_1,\ldots,\sigma_n$.
Since $G$ is compact, $\sigma$ is proper and separates orbits of $G$,
it thus induces a homeomorphism between $V/G$ and $\sigma(V)$, by the 
following lemma.

\begin{lemma*}
Suppose that $X$ and $Y$ are locally compact, Hausdorff spaces and that 
$f : X \to Y$ is bijective, continuous, and proper. Then $f$ is a 
homeomorphism.
\end{lemma*} 

\proof
(E.g.\!\! \cite{bredontopgeo})
By defining $\tilde f(\infty)=\infty$, $f$ extends to a continuous 
map $\tilde f : X \cup \{\infty\} \to Y \cup \{\infty\}$ between the 
one point compactifications, since it is proper. 
If $A \subseteq X$ is closed in $X$, then 
$A \cup \{\infty\}$ is closed in $X \cup \{\infty\}$ and hence compact. 
Then, $\tilde f(A \cup \{\infty\})$ is compact and hence closed in 
$Y \cup \{\infty\}$. Consequently, 
$f(A) = \tilde f(A \cup \{\infty\}) \cap Y$ is closed in $Y$.
\endproof

\subsection{Description of $\sigma(V)$} \label{sigmaV}

Let $\< \quad \mid \quad \>$ denote also the $G$-invariant dual 
inner product on $V^*$. The differentials $d \sigma_i : V \to V^*$ are 
$G$-equivariant, and the polynomials 
$v \mapsto \< d \sigma_i(v) \mid d \sigma_j(v) \>$ are in $\mathbb R[V]^G$ 
and are entries of an $n \times n$ symmetric matrix valued polynomial
\[
B(v) := 
\left(
\begin{array}{ccc}
\< d \sigma_1(v) \mid d \sigma_1(v) \> & \cdots & 
\< d \sigma_1(v) \mid d \sigma_n(v) \> \\
\vdots & \ddots & \vdots \\
\< d \sigma_n(v) \mid d \sigma_1(v) \> & \cdots & 
\< d \sigma_n(v) \mid d \sigma_n(v) \> 
\end{array}
\right). 
\]
There is a unique matrix valued polynomial $\tilde B$ on $Z$ such that 
$B = \tilde B \circ \sigma$. The following theorem is due to Procesi and Schwarz 
\cite{procesischwarz}.

\begin{theorem*}
$\sigma(V) = \{z \in Z : \tilde B(z) ~\mbox{positive semidefinite}\}$.
\end{theorem*}

This theorem provides finitely many equations and inequalities 
describing $\si(V)$. Changing the choice of generators may change 
the equations and inequalities, but not the set they describe.

For each $1 \le i_1 < \cdots < i_s \le n$ and 
$1 \le j_1 < \cdots < j_s \le n$ ($s \le n$) consider the matrix with entries 
$\langle d \sigma_{i_p} \mid d \sigma_{j_q} \rangle$ for $1 \le p,q \le s$. 
Denote its determinant by $\Delta_{i_1,\ldots,i_s}^{j_1,\ldots,j_s}$. 
Then, $\Delta_{i_1,\ldots,i_s}^{j_1,\ldots,j_s}$ is a $G$-invariant polynomial 
on $V$, and thus there is a unique polynomial 
$\tilde{\Delta}_{i_1,\ldots,i_s}^{j_1,\ldots,j_s}$ on $Z$ such that 
$\Delta_{i_1,\ldots,i_s}^{j_1,\ldots,j_s} 
= \tilde{\Delta}_{i_1,\ldots,i_s}^{j_1,\ldots,j_s} \circ \sigma$. 

\subsection{The problem of lifting curves} \label{liftingproblem}
Let $c : \mathbb{R} \to V/G = \sigma(V) \subseteq \mathbb{R}^n$ be a smooth curve 
in the orbit space; smooth as curve in $\mathbb{R}^n$.
A curve $\bar{c} : \mathbb{R} \to V$ is called lift of $c$ to $V$, if $c = \sigma \circ \bar{c}$ 
holds. 
\emph{The problem of lifting smooth curves over invariants is independent of the 
choice of a system of homogeneous generators of $\mathbb{R}[V]^G$} 
in the following sense: Suppose $\sigma_1,\ldots,\sigma_n$ and 
$\tau_1,\ldots,\tau_m$ both generate $\mathbb{R}[V]^G$. 
Then for all $i$ and $j$ we have 
$\sigma_i = p_i(\tau_1,\ldots,\tau_m)$ and 
$\tau_j = q_j(\sigma_1,\ldots,\sigma_n)$ for polynomials $p_i$ and $q_j$. 
If $c^{\sigma} = (c_1,\ldots,c_n)$ is a curve in $\sigma(V)$, then 
$c^{\tau} = (q_1(c^{\sigma}),\ldots,q_m(c^\sigma))$ defines a curve in 
$\tau(V)$ of the same regularity. Any lift $\bar{c}$ to $V$ of the curve 
$c^{\sigma}$, i.e., $c^{\sigma} = \sigma \circ \bar{c}$, is a lift of 
$c^{\tau}$ as well (and conversely):
\[c^{\tau} = (q_1(c^{\sigma}),\ldots,q_m(c^\sigma))
=(q_1(\sigma(\bar{c})),\ldots,q_m(\sigma(\bar{c})))
=(\tau_1(\bar{c}),\ldots,\tau_m(\bar{c})) 
=\tau \circ \bar{c}.
\]

\subsection{Stratification of the orbit space} \label{stratification} 

Let $H = G_v$ be the isotropy group of $v \in V$ and $(H)$ the conjugacy class 
of $H$ in $G$ which is called the type of an orbit $G.v$. The union $V_{(H)}$ 
of orbits of type $(H)$ is called an orbit type submanifold of the 
representation $\rho$ and $V_{(H)}/G$ is called an orbit type submanifold of 
the orbit space $V/G$. 
The collection of connected components of the manifolds $\{V_{(H)}/G\}$ 
forms a stratification of $V/G$ called orbit type stratification, 
see \cite{Pflaum}, \cite{schwarz1}.  
The semialgebraic subset $\sigma(V) \subseteq \mathbb R^n$ is naturally 
Whitney stratified (\cite{Loja}). 
The homeomorphism of $V/G$ and $\sigma(V)$ induced by 
$\sigma$ provides an isomorphism between the orbit type stratification of $V/G$ 
and the primary Whitney stratification of $\sigma(V)$, 
see \cite{bierstone}. 
These facts are essentially consequences of the slice theorem, see 
e.g.\! \cite{schwarz1}.

The inclusion relation on the set of subgroups of $G$ induces a partial ordering on the family of conjugacy classes. 
There is a unique minimum orbit type, the principal orbit type, corresponding 
to the open and dense submanifold $V_\reg$ (respectively $V_\reg/G$) 
consisting of regular points, i.e., points where the isotropy 
representation is trivial. The points in the complement $V_\sing$ 
(respectively $V_\sing/G$) are called singular.

\begin{theorem*} \cite{procesischwarz}
Let $\tilde B$ be as in \ref{sigmaV}. The $k$-dimensional primary 
strata of $\sigma(V)$ are the connected components of the set 
$\{z \in \sigma(V) : \on{rank} \tilde B(z) = k\}$.
\end{theorem*}

\subsection{Smooth lifts} \label{smoothlifts}

Let us recall some results from \cite{rep-lift}.

Let $s \in \mathbb{N}_0$. 
Denote by $A_s$ the union of all strata $X$ of the orbit 
space $V/G$ with $\dim X \le s$, and by $I_s$ the ideal of 
$\mathbb{R}[Z] = \mathbb{R}[V]^G$ consisting of all polynomials 
vanishing on $A_{s-1}$.
Let $c : \mathbb{R} \rightarrow V/G = \sigma(V) \subseteq \mathbb{R}^n$ be
a smooth curve, $t \in \mathbb{R}$, and $s = s(c,t)$ a minimal integer 
such that, for a neighborhood $J$ of $t$ in $\mathbb{R}$, 
we have $c(J) \subseteq A_s$. 
The curve $c$ is called \emph{normally nonflat at $t$} 
if there is $f \in I_s$ such that $f \circ c$ is nonflat at $t$,
i.e., the Taylor series of $f \circ c$ at $t$ is not identically zero.
A smooth curve $c : \mathbb R \to \sigma(V) \subseteq \mathbb R^n$ is 
called \emph{generic}, if $c$ is normally nonflat at $t$ for each 
$t \in \mathbb R$.
 
It is easy to see, that $c$ is normally nonflat at $t \in \mathbb R$
if there is some integer $1 \le r \le n$ such that:
\begin{enumerate}
\item The functions 
$\tilde{\Delta}_{i_1,\ldots,i_k}^{j_1,\ldots,j_k} \circ c$ 
vanish in a neighborhood of $t$ whenever $k > r$.
\item There exists a minor  $\tilde{\Delta}_{i_1,\ldots,i_r}^{j_1,\ldots,j_r}$ 
such that $\tilde{\Delta}_{i_1,\ldots,i_r}^{j_1,\ldots,j_r} \circ c$ is 
nonflat at $t$.
\end{enumerate}

\begin{theorem*} 
Let $c : \mathbb{R} \rightarrow \sigma(V) \subseteq \mathbb{R}^n$ be a
smooth curve which is normally nonflat at $t \in \mathbb R$. 
Then there exists a smooth lift $\bar c$ in $V$ of $c$, 
locally near $t$.
If $c$ is generic then there exists a global smooth lift 
$\bar c$ of $c$.
\end{theorem*}

\subsection{Smooth roots} \label{smoothroots}

In the special case that the symmetric group $\on{S}_n$ is acting on 
$\mathbb R^n$ by permuting the coordinates there is the following 
interpretation of the described lifting problem. 
As generators of $\mathbb{R}[\mathbb R^n]^{\on{S}_n}$ we may take the 
elementary symmetric functions
\[
E_j(x) = 
\sum_{1 \le i_1 < \cdots < i_j \le n} x_{i_1} \cdots x_{i_j}
\qquad (1 \le j \le n),
\]
which constitute the coefficients $a_j$ of a monic polynomial 
\[
P(x)=x^n-a_1 x^{n-1}+ \cdots + (-1)a_n
\] 
with roots $x_1,\ldots,x_n$ via Vieta's formulas.
Then a curve in the orbit space $\mathbb R^n/\on{S}_n = E(\mathbb R^n)$
corresponds to a curve $P(t)$ of monic polynomials of degree $n$ with 
only real roots (such polynomials are called \emph{hyperbolic}), 
and a lift of $P(t)$ may be interpreted as a 
parameterization of the roots of $P(t)$. 

The first $n$ Newton polynomials 
\[
N_i(x_1,\ldots,x_n) = \sum_{j=1}^n x_j^i
\] 
which are related to the elementary symmetric functions by
\begin{equation} \label{rec}
N_k - N_{k-1} E_1 + N_{k-2} E_2 + \cdots + (-1)^{k-1} N_1 E_{k-1} 
+ (-1)^k k E_k = 0 \quad (k \ge 1)
\end{equation}
constitute a different system of generators of $\mathbb R[\mathbb R^n]^{\on{S}_n}$.
For convenience we shall switch from elementary symmetric functions to Newton 
polynomials and conversely, if it seems appropriate.

Let us choose $\frac{1}{j} N_j$, $1 \le j \le n$, as generators of 
$\R[\R^n]^{\on{S}_n}$ and put 
$\De_k := \De_{1,\ldots,k}^{1,\ldots,k}$ and 
$\tilde \De_k := \tilde \De_{1,\ldots,k}^{1,\ldots,k}$. 
Then (\cite{roots})
\begin{equation} \label{De}
\De_k(x)=\sum_{i_1<\cdots<i_k} 
(x_{i_1}-x_{i_2})^2 \cdots (x_{i_1}-x_{i_k})^2 \cdots (x_{i_{k-1}}-x_{i_k})^2. 
\end{equation}  

\begin{theorem*}
\cite{roots}
Consider a smooth curve $P(t)$, $t \in \mathbb R$, of monic hyperbolic 
polynomials of fixed degree $n$. Let one of the following two 
equivalent conditions be satisfied:
\begin{enumerate}
\item If two of the increasingly ordered continuous roots meet of 
infinite order at $t_0$ then their germs at $t_0$ are equal.
\item Let $k$ be maximal with the property that the germ at $t_0$ of
$\tilde \Delta_k(P)$ is not $0$.
Then $\tilde \Delta_k(P)$ is not infinitely flat at $t_0$.
\end{enumerate}
Then $P(t)$ is smoothly solvable near $t = t_0$. 
If $(1)$ or $(2)$ are satisfied for any $t_0 \in \mathbb R$, then 
the roots of $P$ may be chosen smoothly globally, and any two 
choices differ by a permutation.
\end{theorem*}

\begin{lemma*}
Condition $(1)$ (and thus condition $(2)$) in the above theorem is 
satisfied if and only if $P$ is normally nonflat at $t_0$ as curve 
in $E(\R^n) = \R^n/\on{S}_n$.
\end{lemma*}

\proof
Let $P$ be normally nonflat at $t_0$.
Let $s$ be a minimal integer such that $P(t)$ lies in $A_s$ for $t$ near $t_0$ 
and let $f \in I_s$ be such that $f \circ P$ is not infinitely flat at $t_0$. 
Denote by $\bar I_s$ the ideal in $\R[\R^n]$ defining the closed subset 
$\pi^{-1}(A_{s-1}) \subseteq \R^n$, where $\pi: \R^n \to \R^n/\on{S}_n$ is the 
quotient projection. 
It is easy to see that the polynomials 
\[
f_{i_1\dots i_s}=(x_{i_1}-x_{i_2})\cdots(x_{i_1}-x_{i_s})\cdots(x_{i_{s-1}}-x_{i_s}),
\] 
where $1\le i_1<\dots <i_s\le n$, generate $\bar I_s$. 
So there exist polynomials 
$Q_{i_1\dots i_s}\in\R[\R^n]$ such that 
\[
f \circ \pi = \sum_{i_1<\dots<i_s} Q_{i_1\dots i_s} f_{i_1\dots i_s}.
\]
Denote by $\bar P(t)$ the lift of $P(t)$ given by the increasingly ordered continuous 
roots $x_1(t),\dots,x_n(t)$ of the polynomial $P(t)$. Then we have 
\[
f \circ P(t)=
\sum_{i_1<\dots<i_s}Q_{i_1\dots i_s} \circ \bar P(t)\cdot f_{i_1\dots i_s} \circ \bar P(t).
\]
Since $f \circ P$ is not infinitely flat at $t_0$, 
at least one of the summands in this sum is not infinitely 
flat at $t_0$ and thus there is a polynomial $f_{i_1\dots i_s}$ such that 
$f_{i_1\dots i_s} \circ \bar P$ is not infinitely flat at $t_0$. 
By assumption, among the roots $x_1(t),\dots, x_n(t)$ there are precisely $s$ 
distinct for $t$ near $t_0$. Hence the germs at $t_0$ of the roots 
$x_{i_1}(t),\dots,x_{i_s}(t)$ are distinct, 
and no two of them meet of infinite order at $t_0$. 
Therefore, condition $(1)$ in the above theorem is satisfied. 

The other direction is evident by \eqref{De}.
\endproof

\section{Lifting smooth curves in spaces of hyperbolic polynomials} 
\label{orbemb}

\subsection{The problem} \label{problem}

Let us denote by $\Hyp^n$ the space of hyperbolic polynomials of 
degree $n$ 
\[
P(x) = x^n + \sum_{j=1}^n (-1)^j a_j x^{n-j}.
\]
We may naturally view $\Hyp^n$ as a semialgebraic subset of 
$\mathbb R^n$ by identifying $P$ with $(a_1,\ldots,a_n)$. 
We have $\Hyp^n = E(\mathbb R^n) = \mathbb R^n/\on{S}_n$, and, 
by means of \ref{sigmaV}, we may calculate explicitly a set of 
inequalities defining $\Hyp^n$ (no equalities since 
the ring $\mathbb R[\mathbb R^n]^{\on{S}_n}$ is polynomial).  

Suppose $X$ is a semialgebraic subset of $\Hyp^n$. 
Let $c : \mathbb R \to X$ be a smooth curve in $X$; 
smooth as curve in $\mathbb R^n$. We may view $c$ as a curve in 
$\Hyp^n$, i.e., as a smooth curve of monic hyperbolic polynomials 
of degree $n$. 
In \ref{smoothroots} sufficient conditions for the existence 
of a smooth lift $\bar c$ to $\mathbb R^n$, i.e., a smooth 
parameterization of its roots, are presented.
It is evident that a smooth curve $c$ in $X$ in order to be liftable 
smoothly over $E$ to $E^{-1}(X)$ must in general fulfill weaker 
genericity conditions. Our purpose is to investigate 
that phenomenon. 

\subsection{Orbit spaces embedded in spaces of hyperbolic polynomials}
\label{chern}

We recall a construction due to L. Smith and R.E. Strong 
\cite{smithstrong} (see also \cite{barbanconrais}) related to E. Noether's 
\cite{noether} proof of Hilbert's finiteness theorem as recounted by 
H. Weyl \cite{weyl}.

Let $\rho : G \to \on{GL}(V)$ be a representation of a finite group $G$ in a
finite dimensional vector space $V$. 
Consider its induced representation in the dual $V^*$. 
For an orbit $B \subseteq V^*$ set 
\[
\phi_B(X) = \prod_{b \in B} (X+b)
\]
which we regard as an element of the ring $\mathbb R[V][X]$, with 
$X$ a new variable.
The polynomial $\phi_B(X)$ is called the \emph{orbit polynomial} of $B$. 
Evidently, $\phi_B \in \mathbb R[V]^G[X]$.
If $|B|$ denotes the cardinality of the orbit $B$, we may expand 
$\phi_B(X)$ to a polynomial of degree $|B|$ in $X$,
\[
\phi_B(X) = \sum_{i+j=|B|} C_i(B) X^j,
\]
defining classes $C_i(B) \in \mathbb R[V]^G$ called the 
\emph{orbit Chern classes} of $B$.

\begin{theorem*}
[L. Smith and R.E. Strong \cite{smithstrong}]
Let $\rho : G \hookrightarrow \on{GL}(V)$ be a faithful representation of a 
finite group $G$. Then there exist orbits $B_1,\ldots,B_l \subseteq V^*$ 
such that the associated orbit Chern classes $C_i(B_j)$,
$1 \le i \le |B_j|$, $1 \le j \le l$, generate $\mathbb R[V]^G$.
\end{theorem*}

The field of real numbers may be replaced by any field of either 
characteristic zero or characteristic larger than the order of $G$. 
For our purpose the reals will suffice.

The Chern classes of the orbit are exactly the elementary symmetric 
functions in the elements of the orbit. 
If $B \subseteq V^*$ is an orbit and $V_B^*$ is a vector space with 
basis identified with the elements of $B$, 
then there is a natural map $V_B^* \to V^*$ given by the identification.
This map induces a map $\mathbb R[V_B]^{\on{S}_{|B|}} \to \mathbb R[V]^G$ 
which sends the $k$-th elementary symmetric function to the $k$-th 
orbit Chern class of $B$.    

In this notation the above theorem says that there exist orbits 
$B_1,\ldots,B_l \subseteq V^*$ such that the induced map 
\[
\bigotimes_{i=1}^l \mathbb R[V_{B_i}]^{\on{S}_{|B_i|}} \longrightarrow 
\mathbb R[V]^G
\]
is surjective.

The orbit Chern classes $C_i(B)$ of an orbit $B$, viewed as invariant 
polynomials on $V$, define a $G$-invariant map 
\[
C(B) = (C_1(B),\ldots,C_{|B|}(B)) : V \longrightarrow \mathbb R^{|B|}
\]
whose image $C(B)(V)$ is a semialgebraic subset of the space $\Hyp^{|B|}$ 
of hyperbolic polynomials of degree $|B|$.

According to \ref{setting} and the above theorem, for any faithful
representation $\rho : G \hookrightarrow \on{GL}(V)$ of a finite group $G$ 
there exist orbits $B_1,\ldots,B_l \subseteq V^*$ such that the map 
\[
C(\rho) = (C(B_1),\ldots,C(B_l)) : V \longrightarrow 
\Hyp^{|B_1|} \times \cdots \times \Hyp^{|B_l|} \subseteq 
\mathbb R^{|B_1|+\cdots+|B_l|}
\]
induces a homeomorphism between the orbit space $V/G$ and the 
image $C(\rho)(V)$ which is a semialgebraic subset of 
$\Hyp^{|B_1|} \times \cdots \times \Hyp^{|B_l|}$. 
By increasing the number of orbits $B_i$ if necessary, we may 
assume that each irreducible subspace of $V$ contributes at least one 
orbit $B_i$. Then, the linear forms $b \in B_1 \cup \cdots \cup B_l$ 
induce an injective inclusion 
$V \hookrightarrow \mathbb R^{|B_1|+\cdots+|B_l|}$.

Let $c : \mathbb R \to C(\rho)(V)$ be a smooth curve. 
Then $c = (c_1,\ldots,c_l)$ where each 
$c_i : \mathbb R \to C(B_i)(V)$ is smooth. 
Since $C(B_i)(V) \subseteq \Hyp^{|B_i|}$ we may view $c_i$ as a 
curve in $\Hyp^{|B_i|}$. 
If there exist smooth lifts $\bar c_i : \mathbb R \to \mathbb R^{|B_i|}$  
with respect to the representations $\on{S}_{|B_i|} : \mathbb R^{|B_i|}$, 
then 
$\bar c = (\bar c_1,\ldots,\bar c_l) : \mathbb R \to \mathbb R^{|B_1|+\cdots+|B_l|}$ is a smooth lift with respect to 
$\on{S}_{|B_1|} \times \cdots \times \on{S}_{|B_l|} 
: \mathbb R^{|B_1|+\cdots+|B_l|}$.
Consequently, it suffices to study the case when there is given a smooth 
curve in a semialgebraic subset of some $\Hyp^n$. 
That is exactly the problem introduced in \ref{problem}.   

Suppose $\tilde c : \mathbb R \to V$ is a smooth lift of $c$ with 
respect to $\rho$. 
Then, there exists a smooth lift 
$\bar c : \mathbb R \to \mathbb R^{|B_1|+\cdots+|B_l|}$ of $c$ with 
respect to the representation of 
$\on{S}_{|B_1|} \times \cdots \times \on{S}_{|B_l|}$ on 
$\R^{|B_1|+\cdots+|B_l|}$, namely 
\[
\xymatrix{
& V \ar@{^{(}->}[r] \ar[d] & \mathbb R^{|B_1|+\cdots+|B_l|} \ar[d] \\ 
\mathbb R \ar[r]_(.35){c} \ar[ur]^{\tilde c} & C(\rho)(V) \ar@{^{(}->}[r] & 
\Hyp^{|B_1|} \times \cdots \times \Hyp^{|B_l|}
}
\]
It follows, by \ref{smoothlifts}, that conditions which guarantee that 
$c$ is generic as curve in the orbit space $V/G$ suffice to imply 
the existence of a smooth lift of $c$ with respect to
$\on{S}_{|B_1|} \times \cdots \times \on{S}_{|B_l|} 
: \mathbb R^{|B_1|+\cdots+|B_l|}$. 

We have seen that the above construction provides a class of 
semialgebraic subsets of spaces of hyperbolic polynomials, 
namely orbit spaces of faithful finite group representations, for which 
we are able to apply the strategy described in \ref{problem}, 
thanks to the results of \ref{smoothlifts}.

In the remaining sections we shall change the point of view. 
Assume we are given a curve of hyperbolic polynomials with certain 
symmetries. We will investigate whether we can weaken the 
conditions in \ref{smoothroots} which guarantee the existence of 
smooth parameterizations of the roots. 
This will be performed in section \ref{polysym}. 
The following section provides the necessary preparation.

\section{Orbit type and ambient stratification} \label{orbamb}

Suppose $U$ is a linear subspace of $\mathbb R^n$. 
Let the symmetric 
group $\on{S}_n$ act on $\mathbb R^n$ by permuting the coordinates and 
endow $U$ with the induced effective action of 
\[
W = W(U) := N(U)/Z(U),
\]
where $N(U) := \{\tau \in \on{S}_n : \tau.U=U\}$ and 
$Z(U) := \{\tau \in \on{S}_n : \tau.x=x ~\mbox{for all}~ x \in U\}$. 
Then $U$ carries two natural stratifications: 
the orbit type stratification with respect to the $W$-action and 
the restriction to $U$ of the orbit type stratification of $\mathbb R^n$ 
with respect to the $\on{S}_n$-action. It is easily seen that the latter 
indeed provides a Whitney stratification of $U$. 
Let us denote it as the \emph{ambient stratification of $U$}. 

\begin{proposition} \label{stratcomp}
Let $U$ be a linear subspace in $\R^n$ endowed with the induced 
action by $W=W(U)$. Then for the ambient and orbit type stratification of 
$U$ we have: 
\begin{enumerate}
\item Each ambient stratum is contained in a unique orbit type stratum.
\item Each orbit type stratum contains at least one ambient stratum of the 
same dimension and is the union of all contained ambient strata.
\end{enumerate} 
\end{proposition}

\proof 
To $(1)$: 
Let $S$ be an ambient stratum, i.e.,  
$S$ is a component of $\on{S}_n.\R^n_H \cap U$, 
where $H = (\on{S}_n)_x$ for a $x \in U$ and 
$\R^n_H=\{y \in \R^n : (\on{S}_n)_y = H\}$. 
Since $\on{S}_n$ is finite and the manifolds $\tau.\R^n_H$ for $\tau \in \on{S}_n$ 
either coincide or are pairwise disjoint, the components of $\on{S}_n.\R^n_H$ 
are open subsets of $\tau.\R^n_H$ for $\tau \in \on{S}_n$. 
Thus, we may assume that $S$ is a component of $\R^n_H \cap U$.

Denote by $\pi$ the quotient projection $N(U)\to N(U)/Z(U)=W$. 
For any $u \in U$ we have $W_u=\pi(N(U)\cap (\on{S}_n)_u)$ and thus 
$\R^n_H \cap U \subseteq \{u \in U : W_u=W_x\}$. 
By definition and a similar argument as above, 
the components of the subset $\{u\in U : W_u=W_x\}$  
are orbit type strata of $U$. 
So the ambient stratum $S$ is contained in a unique isotropy type stratum 
$R_S$. 
  
To $(2)$: Let $R$ be an orbit type stratum and let $\mathfrak S$ be the set 
of all ambient strata $S$ such that $R_{S} = R$, where $R_{S}$ 
is the unique orbit type stratum from $(1)$. 
Clearly, $R = \bigcup \mathfrak S$ 
and for each $S \in \mathfrak S$ we have $\dim S \le \dim R$. 
Since the set $\mathfrak S$ is finite, 
there is a stratum $S \in \mathfrak S$ such that $\dim S = \dim R$. 
\endproof

\begin{remarks}
$(1)$ It is easy to see that proposition \ref{stratcomp} is true if one 
replaces the $\on{S}_n$-module $\R^n$ by any finite dimensional $G$-module $V$, 
where $G$ is a finite group.

$(2)$ Proposition \ref{stratcomp} implies that the orbit type 
stratification of $U$ is coarser than its ambient stratification. 
That means, following \cite{Pflaum}, that for each ambient stratum $S$ 
there exists an orbit type stratum $R_S$ such that $S \subseteq R_S$, 
$\id|_S : S \to R_S$ is smooth, and for all $S \subseteq \overline{S'}$ 
we have $R_S \subseteq \overline{R_{S'}}$. It remains to check the last
condition:  
Assume that $S \subseteq \overline{S'}$.
Since $S \subseteq R_{S}$ and 
$S \subseteq \overline{S'} \subseteq \overline{R_{S'}}$, 
we obtain $R_{S} \cap \overline{R_{S'}} \ne \emptyset$, and, 
by the frontier condition,
$R_{S} \subseteq \overline{R_{S'}}$.
\end{remarks}

Assume that the restrictions $E_i|_U$, $1 \le i \le n$,  
generate the algebra $\R[U]^W$.
It follows that $E|_U = (E_1|_U,\ldots,E_n|_U)$ induces a 
homeomorphism between $U/W$ and the semialgebraic subset $E(U)$ of 
$\R^n/\on{S}_n = E(\R^n) = \Hyp^n$, by \ref{setting}. 
It is well-known that $U_{(H)} \to U_{(H)}/W$, where $H =W_u$ for some $u \in U$, 
is a Riemannian submersion. Since $W$ is finite, it is even a local 
diffeomorphism. By proposition \ref{stratcomp}, this implies that 
for any ambient stratum $S$ in $U$ the image $E(S)$ is a smooth manifold. 
The collection $\mathcal T = \{E(S) : S ~\text{ambient stratum in}~ U\}$
obviously coincides with the collection obtained by restricting to $E(U)$ the 
orbit type stratification of $\R^n/\on{S}_n = E(\R^n) = \Hyp^n$.
It is easily verified that the frontier condition for the orbit type 
stratification of $\R^n/\on{S}_n = E(\R^n) = \Hyp^n$ implies the frontier
condition for $\mathcal T$. 
Consequently, $\mathcal T$ provides a stratification of $E(U)$. 
Let us denote this stratification as the 
\emph{ambient stratification of $E(U)$}. 

Consider a smooth curve $c : \R \to E(U) = U/W$ in the sense of \ref{liftingproblem}.
It may then be also viewed as a smooth curve in $\R^n/\on{S}_n = E(\R^n) = \Hyp^n$. 
Thus it makes sense to speak about the normal nonflatness 
of $c$ at some point $t_0$ with respect to the orbit type stratification 
of $U/W$ on the one hand and with respect to the orbit type 
stratification of $\R^n/\on{S}_n$ on the other hand.
To shorten notation we shall say that $c$ is normally nonflat at $t_0$ 
with respect to the ambient stratification of $U/W$ iff it is 
normally nonflat at $t_0$ with respect to the orbit type stratification of 
$\R^n/\on{S}_n$.

\begin{proposition} \label{nnamborb}
Let $U$ be a linear subspace in $\R^n$ endowed with the induced 
action by $W=W(U)$ and assume that 
the restrictions $E_i|_U$, $1 \le i \le n$, generate $\R[U]^W$. 
Consider a smooth curve $c : \R \to E(U) = U/W$.  
If $c$ is normally nonflat at $t_0$ with respect to the ambient 
stratification of $U/W$, then it is normally nonflat at $t_0$ with respect 
to the orbit type stratification of $U/W$.
\end{proposition}

\proof 
The set of reflection hyperplanes $H$ of the reflection group $\on{S}_n$ is 
in bijective correspondence with the set of linear functionals $\om_H$ 
on $\R^n$ of the form $x_j-x_i$ for $1 \le i < j \le n$, namely $H$ is 
the kernel of $\om_H$. Let us consider the restrictions $\om_H|_U$ to 
$U$.
If $c$ is normally nonflat at $t_0$ with respect to the ambient 
stratification, then, by lemma \ref{smoothroots}, any two of the increasingly 
ordered continuous roots of the polynomial $c(t) \in E(U) \subseteq \Hyp^n$ 
either coincide identically near $t_0$ or do not meet at $t_0$ of infinite 
order. Then for the continuous lift $\bar c$ of $c$ defined by such a choice 
of roots any function $\om_H \circ \bar c$ either vanishes identically near 
$t_0$ or does not vanish at $t_0$ of infinite order. 

Let $s$ be a minimal integer such that $c(t)$ lies in 
$A_{s,\orb}$ for $t$ near $t_0$, where $A_{s,\orb}$ is the union of all 
orbit type strata of $U/W$ of dimension $\le s$.

Denote by $\pi_U$ the projection $U\to U/W$. 
Let $R$ be an orbit type stratum contained in $\pi_U^{-1}(A_{s-1,\orb})$ and 
let $S_1,\dots,S_k$ be the ambient strata of the same dimension as $R$ 
contained in $R$ (see proposition \ref{stratcomp}). 
For each $1\le j \le k$ denote by $\mathcal H_j$ the set of reflection 
hyperplanes for reflections in $\on{S}_n$ fixing $S_j$ pointwise. 
Let $\Om_j$ be the set of linear functionals $\om_H|_U$ for $H \in \mathcal H_j$. 
Put $f_{R,j}=\sum_{\om \in \Om_j}\om^2$. By definition the equation 
$f_{R,j}=0$ defines a linear subspace of $U$ in which $S_j$ is an open subset. 
Let $f_R=\prod_{j=1}^kf_{R,j}$. 
Consider the natural action of $W$ on $\R[U]$ and let 
$W.f_R=\{f_R^1,\ldots,f_R^l\}$ be the orbit through $f_R$ with respect 
to this action. 
Define $F_R=f_R^1\cdots f_R^l$. 
By construction $F_R \in \R[U]^W$ and the set $Z_R$ of zeros of $F_R$ viewed 
as a function on $U/W$ is contained in $A_{s-1,\orb}$. 
Moreover, $A_{s-1,\orb}$ is the union 
of the $Z_R$, where $R$ ranges over all orbit type strata (of maximal dimension) 
contained in $\pi_U^{-1}(A_{s-1,\orb})$. 
Thus $F=\prod_R F_R$, where the product is taken over all orbit type strata
(of maximal dimension) $R$ 
contained in $\pi_U^{-1}(A_{s-1,\orb})$, is a regular function on $U/W$ whose 
set of zeros equals $A_{s-1,\orb}$. 
By construction, the function $F \circ c$ is nonflat at $t_0$.
 
This proves the statement.
\endproof

We define $F_{\amb}(c)$ (resp.\! $F_\orb(c)$) to be the set of all 
$t \in \mathbb R$ such that 
$c$ is normally flat at $t$ with respect to the ambient (resp.\! orbit type) 
stratification of $E(U)$. 
It follows that in the situation of proposition \ref{nnamborb} we have 
$F_{\orb}(c) \subseteq F_{\amb}(c)$.

\section{Choosing roots of polynomials with symmetries smoothly} 
\label{polysym}

Consider a smooth curve of hyperbolic polynomials 
\[
P(t)(x) = x^n - a_1(t) x^{n-1} + a_2(t) x^{n-2} - \cdots + (-1)^n a_n(t)
\qquad (t \in \mathbb R).
\]
We are interested in conditions that guarantee 
the existence of a smooth parameterization of the roots of $P$.
Such conditions have been found in \cite{roots}, 
see \ref{smoothroots}. There no additional assumptions on the polynomials
$P(t)$ have been made. 

In this section we are going to improve those results if the set of 
roots $x_1(t),\dots,x_n(t)$ of $P(t)$ has symmetries 
additional to its invariance under permutations.

Let as assume that the additional symmetries of $P(t)$ 
are given by linear relations between the roots of $P(t)$. 
Otherwise put, 
there is a linear subspace $U$ of $\mathbb R^n$ such that 
$(x_1(t),\dots,x_n(t)) \in U$ for all $t \in \mathbb R$.
Then, the curve $P(t)$ lies in the semialgebraic subset $E(U)$ of 
$\Hyp^n = E(\mathbb R^n) = \mathbb R^n/\on{S}_n$,
the space of hyperbolic polynomials of degree $n$.

The linear subspace $U \subseteq \mathbb R^n$ inherits an effective 
action by the group $W=W(U)$.
 
Let us suppose that the restrictions $E_i|_U$, $1 \le i \le n$, 
generate the algebra $\mathbb R[U]^W$.
Then $E|_U = (E_1|_U,\ldots,E_n|_U)$ induces a 
homeomorphism between $U/W$ and the semialgebraic subset $E(U)$ of 
$\Hyp^n$, by \ref{setting}. 

\begin{lemma}
Consider a continuous curve of hyperbolic polynomials 
\[
P(t)(x) = x^n - a_1(t) x^{n-1} + a_2(t) x^{n-2} - \cdots + (-1)^n a_n(t)
\qquad (t \in \mathbb R).
\] 
Let $U$ be some linear subspace of $\mathbb R^n$ and assume that 
the restrictions $E_i|_U$, $1 \le i \le n$,  
generate the algebra $\mathbb R[U]^{W(U)}$.
Then the following two conditions are equivalent: 
\begin{enumerate}
\item There exists a continuous parameterization $x(t)$ of the roots $x_1(t),\ldots,x_n(t)$
of $P(t)$ such that $x(t) \in U$ for all $t \in \mathbb R$. 
\item $P(t) \in E(U)$ for all $t \in \mathbb R$.
\end{enumerate}
\end{lemma}

\proof
The implication $(1) \Rightarrow (2)$ is trivial.
Suppose that $P(t)$ is a continuous curve in $E(U)$. 
By assumption, we may view $P(t)$ as a 
curve in the orbit space $U/W(U) \cong E(U)$. It allows a 
continuous lift $x(t)$ into $U$, by \cite{rep-lift3} or \cite{mont}, 
which constitutes a parameterization of the roots of $P(t)$.   
\endproof

The smooth curve of polynomials $P(t)$ which lies in $E(U)$ may 
be viewed as a smooth curve in the orbit space $U/W$
in the sense of \ref{liftingproblem}. 
A smooth lift of $P(t)$ over the orbit map $E|_U$ to the $W$-module $U$
provides a smooth parameterization of the roots of the polynomials $P(t)$.

By theorem \ref{smoothlifts}, we may conclude:
If $P(t)$ is normally nonflat at $t=t_0$ with respect to the 
orbit type stratification of 
$E(U)$, then $P(t)$ is smoothly solvable near $t=t_0$.

Consider the closed sets $F_{\amb}(P)$ and $F_{\orb}(P)$, as defined in section 
\ref{orbamb}. By proposition \ref{nnamborb}, the set $F_{\orb}(P)$ is 
contained in $F_{\amb}(P)$. 
We have found that that $P(t)$ is smoothly solvable locally near any 
$t_0 \in \R \backslash F_{\orb}(P)$. 
Any two smooth parameterizations of the roots of $P(t)$ near such a $t_0$ 
differ by a constant permutation, see theorem \ref{smoothroots}. 
Thus the local solutions may be glued 
to a smooth solution on $\R \backslash F_{\orb}(P)$.

It follows from a result in \cite{rootsII} (see also \cite{rep-lift4}) 
that any smooth curve of monic hyperbolic polynomials of fixed degree
allows a global twice differentiable parameterization of its roots. 
By the methods used in \cite{rootsII}, it is easy to combine this with 
the result above in order to get the following theorem.

\begin{theorem} \label{thss}
Consider a smooth curve of hyperbolic polynomials 
\[
P(t)(x) = x^n - a_1(t) x^{n-1} + a_2(t) x^{n-2} - \cdots + (-1)^n a_n(t)
\qquad (t \in \R).
\]
Let $U$ be some linear subspace of $\R^n$ such that: 
\begin{enumerate}
\item The restrictions $E_i|_U$, $1 \le i \le n$, 
generate the algebra $\mathbb R[U]^{W(U)}$. 
\item $P(t) \in E(U)$ for all $t \in \R$.
\end{enumerate} 
Then: 
There exists a global twice differentiable parameterization of the roots 
of $P(t)$ on $\R$ which is smooth on 
$\R \backslash F_{\orb}(P)$. \qed
\end{theorem}

\begin{remark}
The orbit type stratification and the ambient stratification of $E(U)$ 
do in general not coincide, whence theorem \ref{thss} provides an actual 
improvement of the statement of theorem \ref{smoothroots}.
In other words, in general we have $F_\orb(P) \subsetneq F_\amb(P)$. 
It may, for instance, happen that $P(0)$ is regular in $E(U) = U/W$ but 
singular in $\Hyp^n = \mathbb R^n/\on{S}_n$ and $P(t)$ is normally flat at 
$t=0$ with respect to the ambient stratification. 
See examples in section \ref{refl}.
\end{remark}

Let us suppose that a linear subspace $U$ of $\mathbb R^n$ is given. 
It is then a purely computational problem to check whether the assumptions 
we have made in the forgoing discussion are satisfied. 
There are algorithms in computational invariant theory 
(e.g. \cite{DerksenKemper}, \cite{sturmfels}) 
which allow to 
decide whether the restrictions $E_i|_U$, $1 \le i \le n$, generate the 
algebra $\mathbb R[U]^{W(U)}$. 
If the answer is yes, theorem \ref{sigmaV} provides an explicit way 
to describe the semialgebraic subset $E(U) \subseteq \Hyp^n$ by a 
finite set of polynomial equations and inequalities. 
So the condition that the curve $P$ lies in $E(U)$ may again be 
check computationally. The orbit type stratification and the ambient 
stratification of $E(U)$ can be determined explicitly using theorem 
\ref{stratification}. 
Then all ingredients are supplied in order to decide whether the curve 
$P(t)$ is normally nonflat at some $t = t_0$ with respect to the one 
or the other stratification of $E(U)$.

Note that there are refined approaches and algorithms for 
computing the orbit space $V/G$ and its orbit type stratification
of a $G$-module $V$ (when identified with the image of its orbit map). 
In \cite{SV03} rational parameterizations of the strata are 
obtained, while \cite{bayer} provides an algorithm yielding a 
description of each stratum in terms of a minimal number of 
polynomial equations and inequalities, if $G$ is finite. 

We shall carry out that procedure explicitly in example \ref{B3}. 

\section{Choosing roots of polynomials with symmetries differentiably}
\label{diff}

Consider a curve of hyperbolic polynomials 
\[
P(t)(x) = x^n - a_1(t) x^{n-1} + a_2(t) x^{n-2} - \cdots + (-1)^n a_n(t)
\qquad (t \in \mathbb R).
\]
Then the following results are known: 
\begin{result} \label{results1} 
We have:
\begin{enumerate}
\item If all $a_i$ are of class $C^n$, then there exists a differentiable 
parameterization of the roots of $P(t)$ with locally bounded derivative,
\cite{bronshtein}, \cite{wakabayashi}.
\item If all $a_i$ are of class $C^{2n}$, then any differentiable 
parameterization of the roots of $P(t)$ is actually $C^1$, 
\cite{rootsII}, \cite{Mandai}. 
\item If all $a_i$ are of class $C^{3n}$, then there exists a twice 
differentiable parameterization of the roots of $P(t)$,
\cite{rootsII}.
\end{enumerate}
\end{result}

In \cite{rep-lift4} we have proved the following generalizations:
\begin{result} \label{results2}
Let $\rho : G \to \on O(V)$ be a finite dimensional representation of a 
finite group $G$. Let $d = d(\rho)$ be the maximum of the degrees 
of a minimal system of homogeneous generators 
$\sigma_1,\ldots,\sigma_m$ of $\mathbb R[V]^G$. 
Write $V = V_1 \oplus \cdots \oplus V_l$ as orthogonal direct sum of 
irreducible subspaces $V_i$. Define 
$k_i := \min \{|G.v| : v \in V_i \backslash \{0\}\}$, $1\le i \le l$, and 
$k := \max \{d(\rho),k_1,\ldots,k_l\}$.
Let $c : \mathbb R \to V/G = \sigma(V) \subseteq \R^m$ be a curve in 
the orbit space. 
Then: 
\begin{enumerate}
\item If $c$ is of class $C^k$, then there exists a differentiable 
lift of $c$ to $V$ with locally bounded derivative. 
\item If $c$ is of class $C^{k+d}$, then any differentiable lift of $c$ 
is actually of class $C^1$.
\item If $c$ is of class $C^{k+2d}$, then there exists a twice 
differentiable lift of $c$ to $V$.   
\end{enumerate} 
\end{result}

Again we may use these facts in order to improve the results for 
curves $P(t)$ of hyperbolic polynomials with symmetries. 

Let $U$ be some linear subspace of $\mathbb R^n$ such that
the restrictions $E_i|_U$, $1 \le i \le n$,  
generate the algebra $\mathbb R[U]^{W(U)}$, 
and $P(t) \in E(U)$ for all $t \in \mathbb R$.
It follows that we may view $P(t)$ as a curve in the orbit space 
$U/W(U) = E(U)$, and any lift of $P(t)$ over the orbit map $E|_U$ to $U$ 
gives a parameterization of the roots of $P(t)$ of the same regularity.

Provided that the integer $k$, associated to the $W(U)$-module $U$ 
as above, is less 
than the degree $n$ of the polynomials in $P(t)$, we are able, 
using \ref{results2}, to lower 
the degree of regularity in the assumptions of the statements in 
\ref{results1}. We shall give examples in section \ref{refl}.

\section{Construction of a class of examples} \label{constr}

We will present a class of examples which our considerations apply to.

Let $G \subseteq \on O(V)$ be a finite group whose action on the vector space 
$V$ is irreducible and effective.
Choose some non-zero orbit $G.v$. Introducing some numbering we can 
write $G.v =\{g_1.v,\ldots,g_n.v\}$, where $|G.v|=n$ and $g_i \in G$. 
We define a mapping $F_{G,v} : V \to \mathbb R^n$ by 
\[
F_{G,v}(x) := (\<g_1.v \mid x\>,\ldots,\<g_n.v \mid x\>). 
\]
Since the linear span of $G.v$ spans $V$, the mapping $F_{G,v}$ is a
linear isomorphism onto its image $F_{G,v}(V) =: U_{G,v}$. 
The linear space $U_{G,v} \subseteq \R^n$ carries the action of $W_{G,v} := W(U_{G,v})$ 
and a natural $G$-action given by transformations from $W_{G,v}$. 
Since the $G$-action is irreducible, so is the $W_{G,v}$-action. 
Hence $U_{G,v} \subseteq \{y \in \mathbb R^n : y_1 + \cdots + y_n=0\}$.  
Irreducibility and effectiveness of the $G$-action 
induce an injection $G \hookrightarrow W_{G,v}$. Thus we may consider $G$ as a 
subgroup of $W_{G,v}$, and in this picture $F_{G,v}$ is $G$-equivariant.  

\begin{remark}
The linear space $U_{G,v}$ always intersects the submanifold of regular points in the $\on{S}_n$-module $\R^n$. Namely: 
For $1 \le i < j \le n$ we define 
$U_{i,j}=\{F_{G,v}(x) : \< g_i.v \mid x \> = \< g_j.v \mid x \>, x \in V\}$.
By definition, $U_{i,j}$ is a linear subspace of $U_{G,v}$ 
and $\bigcup_{i < j} U_{i,j}$ is the set of singular points of the 
$\on{S}_n$-module $\mathbb R^n$ contained in $U_{G,v}$. 
Since, by definition, $g_i.v \ne g_j.v$ for any $i < j$, we have 
$\dim U_{i,j} = n-1$. Thus, $\bigcup_{i < j} U_{i,j} \ne U_{G,v}$, which 
gives the assertion.
\end{remark}

Put $P_{G,v} := E \circ F_{G,v}$. Then $P_{G,v}$ is proper, since 
$E$ and $F_{G,v}$ are proper.

\begin{lemma} \label{G=W}
Suppose that $P_{G,v}$ separates $G$-orbits. 
Then we have $G = W_{G,v}$.
\end{lemma}

\proof
The groups $G$ and $W_{G,v}$ have the same orbits in $U_{G,v}$. 
For: 
Suppose that $\tau \in W_{G,v}$ and $x,y \in V$ such that $F_{G,v}(y) = \tau.F_{G,v}(x)$. Since $P_{G,v}$ 
separates orbits,
it follows that there exists some $g \in G$ such that $y = g.x$, 
whence $g.F_{G,v}(x) = \tau.F_{G,v}(x)$. 

Now choose $x \in V$ such that $F_{G,v}(x)$ is a regular 
point of the $W_{G,v}$-module $U_{G,v}$. 
The regular points of any effective linear finite group representation are 
precisely those with trivial isotropy groups. 
We may conclude that $x$ is a regular point of the $G$-module $V$. 
So $|W_{G,v}| = |W_{G,v}.F_{G,v}(x)|=|G.x|=|G|$, and thus $G = W_{G,v}$. 
\endproof

If $P_{G,v}$ separates $G$-orbits, then, by lemma \ref{G=W}, the 
$G = W_{G,v}$-modules $V$ and $U_{G,v}$ are equivalent. In particular, 
it follows that the restriction $E|_{U_{G,v}}$ separates 
$W_{G,v}$-orbits,
$F_{G,v}$ induces a homeomorphism between $V/G$ and 
$U_{G,v}/W_{\rho,v}$,
and $F_{G,v}^* : \mathbb R[U_{G,v}]^{W_{G,v}} 
\to \mathbb R[V]^G$ is an algebra isomorphism.

\begin{proposition}
The following conditions are equivalent:
\begin{enumerate}
\item $P_{G,v}$ separates $G$-orbits. 
\item For all $x \in V$ we have 
$F_{G,v}(G.x) = \on{S}_n.F_{G,v}(x) \cap U_{G,v}$.
\item $P_{G,v}$ induces a homeomorphism between $V/G$ and $P_{G,v}(V)$. 
\end{enumerate}
\end{proposition}

\proof
Since $E$ separates $\on{S}_n$-orbits, for each $x \in V$ 
there exists a $z \in \R^n$ such that $E^{-1}(z) = \on{S}_n.F_{G,v}(x)$.
Then the equivalence of $(1)$ and $(2)$ follows from
\[
P_{G,v}^{-1}(z) = F_{G,v}^{-1}(\on{S}_n.F_{G,v}(x)) = 
F_{G,v}^{-1}(\on{S}_n.F_{G,v}(x) \cap U_{G,v}).
\]
%
%
The equivalence of $(1)$ and $(3)$ follows easily from lemma 
\ref{setting}.
\endproof

Note that the introduced construction of $F_{G,v}$ and $P_{G,v}$  
essentially coincides with the construction of orbit Chern classes as 
described in \ref{chern}.

Let us discuss uniqueness of the above construction.
Suppose $G \subseteq \on O(V)$ is a finite group.
Denote by $\aut(G)$ the group of automorphisms of $G$. Let $S$ be the set 
of all reflections belonging to $G$. 
Denote by $\aut(G,S)$ the group of automorphisms of $G$ preserving the set 
$S$. Let $a \in \aut(G,S)$. A diffeomorphism $T : V \to V$ is called 
$a$-equivariant, if $T \circ g = a(g) \circ T$ for any $g \in G$
(cf. \cite{losik}). 

\begin{lemma} \label{aequi}
Suppose $G \subseteq \on O(V)$ is a finite group.
Let $a \in \aut(G,S)$ and let $T : V \to V$ be an $a$-equivariant 
diffeomorphism. Then the isotropy groups of $x$ and $T(x)$ are isomorphic,
for all $x \in V$, $T$ maps orbits onto orbits, and $T$ induces an 
automorphism of the orbit type stratification of $V$.
\end{lemma}

\proof
It is easily seen that $G_{T(x)} = a(G_x)$ and $T(G.x) = G.T(x)$ for all 
$x \in V$. Further, it is evident that $G_x = gHg^{-1}$ if and only if 
$G_{T(x)} = a(g)a(H)a(g)^{-1}$. The statement follows. 
\endproof

Let $c : \R \to V/G = \si(V) \subseteq \R^n$ be a smooth curve and 
$\bar c : \R \to V$ a smooth lift of $c$. The orbit space $V/G$ has a 
smooth structure given by the sheaf 
$C^\infty(V/G) = C^\infty(V)^G$
of smooth $G$-invariant functions on 
$V$. Then $c$ induces a continuous algebra morphism 
$c^* : C^\infty(V/G) \to C^\infty(\R)$ and $\bar c$ induces a continuous 
algebra morphism $\bar c^* : C^\infty(V) \to C^\infty(\R)$ such that 
$c^* = \bar c^* \circ \si^*$. This algebraic lifting problem is equivalent to 
the geometrical one. 
It is evident that to determine $\bar c^*$ it suffices to know the 
images under $\bar c^*$ of some system of global coordinate functions 
$x_1,\ldots,x_m$, where $m = \dim V$. The same is true for $c^*$, and in this 
case we may take the basic invariants $\si_1,\ldots,\si_n$ as global 
coordinates functions, by Schwarz's theorem \cite{schwarz2}. 
If $f : V/G \to V/G$ is a smooth diffeomorphism one can take instead of 
the $\si_i$ the functions $f^*(\si_i)$ with the same result.  
Thus, the problem of smooth lifting is invariant with respect to the group 
of diffeomorphisms of $V/G$. 
Each such diffeomorphism has a smooth lift to $V$ 
which is an $a$-equivariant diffeomorphism, for some $a \in \aut(G,S)$, see 
\cite{losik}. Conversely, any smooth $a$-equivariant diffeomorphism of 
$V$ induces a smooth diffeomorphism of $V/G$, by lemma \ref{aequi}.

Therefore, we may regard two constructions as described above, carried out 
for distinct points $v$ and $w$ in $V$, as equivalent with respect to our 
lifting problem, if there exists a smooth $a$-equivariant diffeomorphism 
$T : V \to V$ with $v = T(w)$, for some $a \in \aut(G,S)$. 

If $T$ is of a particular form, we can even say more.

\begin{proposition} \label{uniqueness}
Suppose $G \subseteq \on O(V)$ is a finite group. 
Let $v,w \in V \backslash \{0\}$. 
If there exists a homothety or an $a$-equivariant linear orthogonal map 
$T : V \to V$, for some $a \in \aut(G,S)$, 
such that $v = T(w)$, then $P_{G,v}(V)$ and $P_{G,w}(V)$ are homeomorphic, 
and $\mathbb R[E_1 \circ F_{G,v},\ldots,E_n \circ F_{G,v}]$ 
and $\mathbb R[E_1 \circ F_{G,w},\ldots,E_n \circ F_{G,w}]$ 
are isomorphic.

Moreover, in both cases, 
the ambient stratifications of $U_{G,v}$ and $U_{G,w}$ 
are isomorphic, i.e., there exists a linear isomorphism 
$U_{G,v} \to U_{G,w}$ mapping strata onto strata.
\end{proposition}

\proof
If $T$ is a homothety, 
then it is equivariant ($a = \id$) and $U_{G,v}=U_{G,w}$. 
If $T$ is $a$-equivariant linear orthogonal, then, by lemma \ref{aequi}, 
the linear subspaces $U_{G,v}$ and $U_{G,w}$ of $\mathbb R^n$ differ 
only by a permutation from $\on{S}_n$. In both cases $P_{G,v}(V)$ and 
$P_{G,w}(V)$ are homeomorphic, and 
$T^* : \mathbb R[E_1 \circ F_{G,v},\ldots,E_n \circ F_{G,v}] \to 
\mathbb R[E_1 \circ F_{G,w},\ldots,E_n \circ F_{G,w}]$ is an algebra 
isomorphism.

The supplement in the lemma follows immediately from the fact 
that $U_{G,v}$ and $U_{G,w}$  
differ only by a permutation of $\on{S}_n$.
\endproof

If $P(t)$ is a smooth curve of hyperbolic polynomials lying in 
$P_{G,v}(V)$ and provided that the polynomials $E_i \circ F_{G,v}$, 
$1 \le i \le n$,
generate $\R[V]^G$, we may apply the results of sections 
\ref{polysym} and \ref{diff}. 

We will investigate the case of finite reflection groups in the next section.

\section{Finite reflection groups} \label{refl}

Suppose $U$ is a linear subspace of $\mathbb R^n$. 
Let the symmetric 
group $\on{S}_n$ act on $\mathbb R^n$ by permuting the coordinates and 
endow $U$ with the induced action of $W = W(U)$. 
We shall assume in this section that $W$ is a finite reflection group.

\begin{remark} \label{hypl}
If $W$ is a finite reflection group, proposition \ref{stratcomp} 
reduces to the following statement: 
{\it Any reflection hyperplane 
of $W$ in $U$ is the intersection with $U$ of some reflection hyperplane 
of $\on{S}_n$ in $\R^n$.}
For: Let $H$ be a reflection hyperplane of $W$ in $U$. 
By proposition \ref{stratcomp}, there exists a ambient stratum $S$ of $U$ 
such that $S \subseteq H$ and $\dim S = \dim H$. 
Obviously, $S \subseteq (\R^n)_\sing \cap U$, and so there are 
reflection hyperplanes $P_1,\ldots,P_l$ of $\on S_n$ in $\R^n$ which 
contain $S$. Since $\dim S = \dim U-1$, there is a $1 \le i \le n$ such 
that $P_i \cap U$ is a hyperplane in $U$. Since $S$ is contained in both 
$H$ and $P_i \cap U$, we have $H = P_i \cap U$.  
\end{remark}

For any finite reflection group $W \subseteq \on O(U)$ we may write 
$U$ as the orthogonal direct sum of $W$-invariant subspaces 
$U_0=U^W,U_1,\ldots,U_l$ such that $W$ is isomorphic to 
$W_0 \times W_1 \times \cdots \times W_l$, where 
$W_i = \{\tau|_{U_i} : \tau \in W\}$. Each $W_i$ ($i \ge 1$) is one 
of the groups (e.g. \cite{humphreys}) 
\begin{gather*}
\on A_m, m\ge 1; \on B_m, m \ge 2; \on D_m, m \ge 4; \on I_2^m, m \ge 5, m \ne 6; \\ 
\on G_2; \on H_3; \on H_4; \on F_4; \on E_6; \on E_7; \on E_8.
\end{gather*}
It follows that $\mathbb R[U]^W \cong \mathbb R[U_1]^{W_1} \otimes 
\cdots \otimes \mathbb R[U_l]^{W_l}$ and 
$U/W \cong U_1/W_1 \times \cdots \times U_l/W_l$. 
A smooth curve $c = (c_1,\ldots,c_l)$ in the orbit space $U/W$ is then 
smoothly liftable to $U$ if and only if, for all $1 \le i \le l$, 
$c_i$ is smoothly liftable to $U_i$. 
Note that the orbit type stratification of $U/W$ coincides with the 
product stratification of the orbit type stratifications $\mathcal Z_i$ 
of the factors $U_i/W_i$, i.e., 
the strata of $U/W$ are $S_1 \times \cdots \times S_l$, 
where $S_i \in \mathcal Z_i$.
Consequently, in order to apply the results of section \ref{polysym} and 
section \ref{diff} we may consider each factor $U_i/W_i$ separately.
So let us assume that $U$ is an irreducible $W$-module.

To this end we have to check whether the restrictions $E_i|_U$, 
$1 \le i \le n$, generate the algebra $\mathbb R[U]^W$.
In practice this is easily accomplishable: 
The unique degrees $d_1,\ldots,d_m$, where $m = \dim U$, of 
the elements in a minimal system of homogeneous generators of 
$\R[U]^W$ are well known.  
It suffices to compute the Jacobian $J$ of the polynomials 
$E_{d_i}|_U$, $1 \le i \le m$. If $J \ne 0 \in \R[U]$ then 
they generate $\R[U]^W$.
Note that a necessary condition for the $E_i|_U$, 
$1 \le i \le n$, to generate $\mathbb R[U]^W$ is that
the degrees $d_1,\ldots,d_m$ must be pairwise 
distinct, see remark \ref{Dm}.

Let us carry out the construction presented in section \ref{constr} for 
finite irreducible reflection groups $G \subseteq \on O(V)$. 
Let $v \in V \backslash \{0\}$. If the polynomials $E_i \circ F_{G,v}$ 
generate the algebra $\R[V]^G$, then $W_{G,v}$ is a finite 
irreducible reflection group as well, by lemma \ref{G=W}. 
 
Fix a system $\Pi$ of simple roots of $G$. 
For any $v$ in the fundamental domain 
$C=\{x \in V : \<x \mid r\> \ge 0 ~\mbox{for all}~ r \in \Pi\}$, 
the isotropy group $G_v$ is generated by the simple reflections it contains
(e.g. \cite{humphreys}). 

\begin{lemma} \label{auto}
Let $G \subseteq \on O(V)$ be a finite reflection group. 
Each automorphism of the corresponding Coxeter diagram $\Gamma(G)$ 
induces an $a$-equivariant orthogonal automorphism of $V$ 
for some $a \in \aut(G,S)$.
\end{lemma}

\proof
(\cite{losik})
Since the vertices in the Coxeter diagram $\Gamma(G)$ 
represent the simple roots of $G$,
an automorphism $\varphi$ of $\Gamma(G)$, defines uniquely an automorphism 
$a_\varphi \in \aut(G,S)$. 
Suppose the simple roots have unit length. 
Since they form a basis for $V$ the automorphism $\varphi$ defines 
naturally an orthogonal automorphism $T_{\varphi}$ of $V$.
It is easily checked that $T_{\varphi}$ is $a_\varphi$-equivariant.
\endproof

\begin{theorem} \label{generate}
Suppose $G \subseteq \on O(V)$ is a finite irreducible reflection group. 
Let $v \in V \backslash \{0\}$ such that the cardinality of $G_v$ is 
maximal. Then: 
The polynomials $E_i \circ F_{G,v}$, 
$1 \le i \le n$, generate $\mathbb R[V]^G$ and $P_{G,v}$ induces a 
homeomorphism between $V/G$ and $P_{G,v}(V)$ if and only if 
$G \ne \on D_m$, $m \ge 4$.
\end{theorem}

\proof 
By proposition \ref{uniqueness} and lemma \ref{auto}
it suffices to check the statement for one 
single $v \ne 0$ with maximal $G_v$. Choosing $e_1+\cdots+e_m-m e_{m+1}$,
$e_1$, and $e_1$ for $\on A_m$, $\on B_m$, and $\on I_2^m$, respectively,  
one obtains the usual systems of basic invariants. The choice $e_1$ for 
$\on D_m$ yields $F_{\on D_m,e_1}=F_{\on B_m,e_1}$, whence the polynomials 
$E_i \circ F_{\on D_m,e_1}$, $1 \le i \le n=2m$, cannot separate $\on D_m$-orbits.
For the remaining irreducible reflection groups the necessary computations 
have been carried out by Mehta \cite{mehta}. 
\endproof

\begin{remark} \label{Dm}
If for $\on D_{m}$ with $m$ odd one chooses $v = e_1 + \cdots + e_m$, then the 
polynomials $E_i \circ F_{\on D_m,v}$, $1 \le i \le n=2^{m-1}$, 
generate $\mathbb R[\mathbb R^m]^{\on D_m}$, since 
the Jacobian of the polynomials  
$N_{i} \circ F_{\on D_m,w}$, $i = 2,4,\cdots,2n-2,n$, is up to a 
constant factor given by $\prod_{i < j} (x_i^2-x_j^2)$. 
If $m (\ge 4)$ is even, this cannot be true since there have to be two 
basic invariants of degree $m/2$.
\end{remark}

The following theorem is a corollary of theorem \ref{generate} 
and theorem \ref{thss}.

\begin{theorem} \label{minorb}
Suppose $G \subseteq \on O(V)$ is a finite irreducible reflection group and 
$G \ne \on D_m$, $m \ge 4$. 
Let $v \in V \backslash \{0\}$ such that the cardinality of $G_v$ is 
maximal. Let 
\[
P(t)(x) = x^n - a_1(t) x^{n-1} + a_2(t) x^{n-2} - \cdots + (-1)^n a_n(t)
\qquad (t \in \mathbb R)
\] 
be a smooth curve of hyperbolic polynomials of degree $n = |G.v|$ 
lying in $P_{G,v}(V)$ for all $t \in \mathbb R$.
Then there exists a global twice differentiable parameterization of the 
roots of $P(t)$ on $\mathbb R$ which is smooth on 
$\mathbb R \backslash F_\orb$. \qed
\end{theorem}

\begin{remark}
It is easy to see that, under the assumption that the cardinality of $G_v$ 
is maximal, the orbit type stratification and the ambient stratification 
of $U_{G,v}$ coincide only for $G=\on A_m,\on B_m,\on I_2^m$. 
In general, 
if $|G_v|$ is not maximal, 
the orbit type stratification of $U_{G,v}$ will be strictly 
coarser than its ambient stratification. 
\end{remark}

It is easy to compute the integer $k$, associated to orthogonal 
representations of finite groups $G$ in \ref{results2}, 
if $G$ is a finite 
irreducible reflection group. See figure \ref{fig1}.

\begin{center}
\begin{figure}[h]
\begin{tabular}{|c||c|c|c|c|c|c|c|c|c|c|c|}
\hline
$G$ & $\on A_m$ & $\on B_m$ & $\on D_m$ & $\on I_2^m$ & 
$\on G_2$ & $\on H_3$ & $\on H_4$ & $\on F_4$ & $\on E_6$ & $\on E_7$ & $\on E_8$ \\
\hline 
$k$ & $m+1$ & $2m$ & $2m$ & $m$ & $6$ & $12$ & $120$ & $24$ & $27$ & $56$ & $240$ \\
\hline
\end{tabular}
\caption{Irreducible Coxeter groups with associated integer $k$.}
\label{fig1}
\end{figure}
\end{center}

In the situation of theorem \ref{minorb} the strategy discussed in section 
\ref{diff} will lead to no improvement, since $k = n$ by definition. 
But, if we choose $v \in V \backslash \{0\}$ such that $|G_v|$ is not 
maximal, then $k < n$ and the methods of section \ref{diff} will 
yield refinements.

In many cases the following theorem provides an improvement of 
\ref{results1}.

\begin{theorem} \label{reflfin}
Suppose $G \subseteq \on O(V)$ is a finite irreducible reflection group.
Choose some $v \in V \backslash \{0\}$. Put $n = |G.v|$ and let $k$ 
be as in figure \ref{fig1}. 
Suppose that the restrictions $E_i|_{U_{G,v}}$, $1 \le i \le n$, 
generate $\R[U_{G,v}]^{W_{G,v}}$.
Let 
\[
P(t)(x) = x^n - a_1(t) x^{n-1} + a_2(t) x^{n-2} - \cdots + (-1)^n a_n(t)
\qquad (t \in \R)
\] 
be a curve of hyperbolic polynomials lying in $P_{G,v}(V)$ for all 
$t \in \R$. 
Then: 
\begin{enumerate}
\item If all $a_i$ are of class $C^k$, then there exists a differentiable 
parameterization of the roots of $P(t)$ with locally bounded derivative.
\item If all $a_i$ are of class $C^{k+d}$, then any differentiable 
parameterization of the roots of $P(t)$ is actually $C^1$. 
\item If all $a_i$ are of class $C^{k+2d}$, then there exists a twice 
differentiable parameterization of the roots of $P(t)$. \qed
\end{enumerate}
\end{theorem}  

\begin{example} \label{B3}
Consider the Coxeter group $\on B_3$ and choose $v = e_1+e_2+e_3$. We find
\begin{multline*}
F_{\on B_3,v}(x) = (x_1+x_2+x_3,-x_1+x_2+x_3,x_1-x_2+x_3,x_1+x_2-x_3, \\
-x_1-x_2+x_3,-x_1+x_2-x_3,x_1-x_2-x_3,-x_1-x_2-x_3)
\end{multline*}
and $U_{\on B_3,v} = \{y \in \mathbb R^8 : y_i + y_j = 0 ~\mbox{for}~ i+j = 9, 
y_1=y_2+y_3+y_4\}$. It is easy to check that 
$N_{2i} \circ F_{\on B_3,v}$, $1 \le i \le 3$, generate 
$\R[\R^3]^{\on B_3}$, by computing their Jacobian. 
It is readily verified that the set of all reflection hyperplanes of 
$W_{\on B_3,v}$ is given by intersecting the following hyperplanes in 
$\mathbb R^8$ with $U_{\on B_3,v}$ (compare with remark \ref{hypl}): 
\[
\{y_1=y_2,y_1=y_3,y_1=y_4,y_1=y_5,y_1=y_6,y_1=y_7,y_2=y_3,y_2=y_4,y_3=y_4\}.
\]
Furthermore, the intersections with $U_{\on B_3,v}$ of the following hyperplanes 
in $\mathbb R^8$,
\[
\{y_1=y_8,y_2=y_7,y_3=y_6,y_4=y_5\},
\] 
are not among the set of reflection hyperplanes of $W_{\on B_3,v}$. 
Therefore, the orbit type stratification of $U_{\on B_3,v}$ is strictly coarser
than its ambient stratification. 

We follow the recipe for computing orbit type and ambient stratification of 
$E(U_{\on B_3,v})=N(U_{\on B_3,v})$ given at the end of section \ref{polysym}. 
We will present only the outcome of the calculations.
Using $N_{2i} \circ F_{\on B_3,v}$, $1 \le i \le 3$, as basic invariants of 
$\R[\R^3]^{\on B_3}$, we find that the symmetric matrix
$\tilde B = (\tilde b_{ij})$ from \ref{sigmaV} has entries
\begin{align*}
\tilde b_{11} &= 32 z_2,~ 
\tilde b_{12} = 64 z_4,~
\tilde b_{13} = 96 z_6,~
\tilde b_{22} = -3 z_2^3+36 z_2 z_4 + 32 z_6, \\
\tilde b_{23} &= \frac{1}{8}(5 z_2^4-108 z_2^2 z_4+192 z_4^2+544 z_2 z_6), \\
\tilde b_{33} &= \frac{1}{64}(27 z_2^5-300 z_2^3 z_4 - 1140 z_2 z_4^2 
+ 1140 z_2^2 z_6 + 7680 z_4 z_6). 
\end{align*}
Put $\tilde \De_{ij} = \det \left(
\begin{array}{cc} 
\tilde b_{ii} & \tilde b_{ij} \\ 
\tilde b_{ji} & \tilde b_{jj} 
\end{array}\right)$ where $i < j$. Then $N(U_{\on B_3,v})$ is the subset in $\R^8$ 
defined by the following relations
\begin{gather*} 
z_2 \ge 0, \tilde \De_{12} \ge 0, \det \tilde B \ge 0 \\
z_1 = z_3 = z_5 = z_7 = 0, \\
384 z_8 = 5 z_2^4 - 72 z_2^2 z_4 + 48 z_4^2 + 256 z_2 z_6.
\end{gather*}
The $3$-dimensional principal orbit type stratum is given by 
\[
R^{(3)}=N(U_{\on B_3,v}) \cap \{z_2 > 0, \tilde \De_{12} > 0, \det \tilde B > 0\}.
\]
Put 
\begin{align*}
\tilde f_1 &= 53 z_2^6 - 840 z_2^4 z_4 + 1680 z_2^2 z_4^2 + 6144 z_4^3 + 
2752 z_2^3 z_6 - 16128 z_2 z_4 z_6 + 9216 z_6^2, \\
\tilde f_2 &= z_2^3 - 12 z_2 z_4 + 32 z_6.
\end{align*} 
There are three $2$-dimensional orbit type strata
\begin{align*}
R_1^{(2)}=N(U_{\on B_3,v}) &\cap \{z_2 > 0, \tilde \De_{12} > 0, \tilde f_1=0\} \\
R_2^{(2)}=N(U_{\on B_3,v}) &\cap \{z_2 > 0, \tilde \De_{12} = 0, 
\tilde \De_{23} > 0, \tilde f_1=0\} \\
R_3^{(2)}=N(U_{\on B_3,v}) &\cap \{z_2 > 0, \tilde \De_{13} > 0, \tilde f_2 = 0\}, 
\end{align*}
the three $1$-dimensional orbit type strata $R_1^{(1)}$, $R_2^{(1)}$, 
$R_3^{(1)}$ are the connected components of 
\begin{align*}
N(U_{\on B_3,v}) &\cap \{z_2 > 0, \tilde \De_{12} = \tilde \De_{13} =  
\tilde \De_{23} = 0\}, 
\end{align*}
and $R^{(0)}=\{0\}$ is the only $0$-dimensional stratum.

The ambient stratification of $N(U_{\on B_3,v})$ is obtained by cutting 
with the surface $\{z_2^2-4z_4=0\}$. 
There are two $3$-dimensional ambient strata 
\[
S_1^{(3)} = R^{(3)} \cap \{z_2^2-4z_4>0\} \quad \mbox{and} \quad
S_2^{(3)} = R^{(3)} \cap \{z_2^2-4z_4<0\},
\]
five $2$-dimensional ambient strata
\begin{align*}
S_1^{(2)} &= R^{(3)} \cap \{z_2^2-4z_4=0\},~
S_2^{(2)} = R_1^{(2)} \cap \{z_2^2-4z_4>0\},\\
S_3^{(2)} &= R_1^{(2)} \cap \{z_2^2-4z_4<0\},~ 
S_4^{(2)} = R_2^{(2)},~ 
S_5^{(2)} = R_3^{(2)},
\end{align*}
four $1$-dimensional ambient strata $S_1^{(1)}=R_1^{(1)}$, 
$S_2^{(1)}=R_2^{(1)}$, $S_3^{(1)}=R_3^{(1)}$, 
$S_4^{(1)} = R_1^{(2)} \cap \{z_2^2-4z_4=0\}$,
and $S^{(0)}=R^{(0)}=\{0\}$ is the only $0$-dimensional ambient stratum. 
See figure \ref{fig2}.

\begin{center}
\begin{figure}[h]
\includegraphics{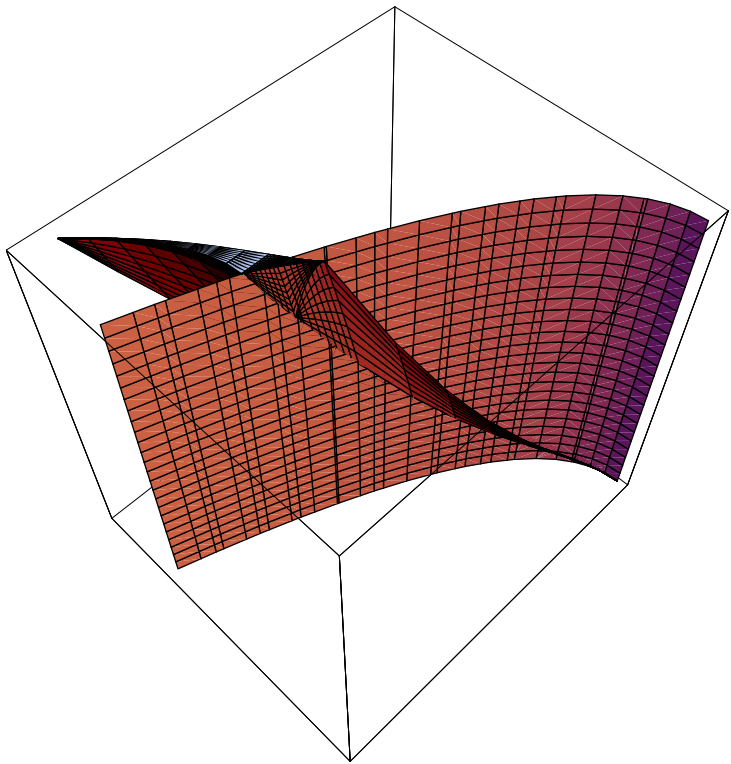}
\caption{The projection of $N(U_{\on B_3,v})$ to the $\{z_2,z_4,z_6\}$-subspace 
and intersection with the surface $\{z_2^2-4z_4=0\}$.}
\label{fig2}
\end{figure}
\end{center}

Let $f$, $g$, $h$ be functions defined in some neighborhood of 
$0 \in \mathbb R$. Suppose that $f$ and $g$ are infinitely flat at $0$ 
and $h(0) = 0$. For $t$ near $0$ consider the curve of polynomials 
$P(t)(x) = x^8 + \sum_{j=1}^8 (-1)^j a_j(t) x^{8-j}$ where 
\begin{gather*}
a_1 = a_3 = a_5 = a_7 = 0, \\
a_2 = -56+f,~ a_4 = 784+g,~ a_6 = -2304+h, \\
1024 a_8 = 16 a_2^4-128 a_2^2 a_4 + 256 a_4^2.
\end{gather*}  
Then, for $t$ near $0$, $P(t)$ is a curve in $N(U_{\on B_3,v})$ 
with $P(0) \in S_1^{(2)}$.
At $t = 0$ it is normally flat with respect to the ambient stratification 
but normally nonflat with respect to the orbit type stratification.

If $f$, $g$ and $h$ are smooth, then $P(t)$ is smoothly solvable 
near $t = 0$, by theorem \ref{thss}.
Note that in this example we have $d=k=6<8=n$ and thus theorem \ref{reflfin} 
provides an actual improvement, too.
\end{example}

The following example shows that $W(U)$ must not necessarily be a finite 
reflection group, even though the $E_i|_U$ generate $\R[U]^{W(U)}$.

\begin{example}
Let $U$ be the subspace of $\R^6$ defined by the following
equations
\[
x_1+x_2+x_3=0,\quad x_4+x_5+x_6=0.
\]
The subgroup $N(U)$ of $\on{S}_6$ is generated by all permutations of 
$x_1,x_2,x_3$, all permutations of $x_4,x_5,x_6$, and the simultaneous 
transpositions of $x_1$ and $x_4$, $x_2$ and $x_5$, $x_3$ and $x_6$. 
The subgroup $Z(U)$ is trivial. 
Thus $W(U)$ is isomorphic to the semidirect product of 
$\on{S}_3 \times \on{S}_3$ and $\on{S}_2$.

One can get the subspace $U$ above as follows.
Consider the point $v=(x,x,x,y,y,y)\in\R^6$, 
where $x,y\ne 0$ and $x\ne y$. 
The isotropy group $H=(\on{S}_6)_v$ of $v$ is evidently isomorphic to 
$\on{S}_3\times \on{S}_3$. Then $U=((\R^6)^H)^\perp$.
The group $H$ is the normal subgroup of $W(U)$ generated by reflections.

First consider the action of $H$ on $U$. 
It is clear that the algebra $\R[U]^H$ is a polynomial
algebra generated by the basic generators 
\begin{gather*}
y_1=x_1^2+x_2^2+x_1x_2,~ z_1=x_1x_2(x_1+x_2),\\ 
y_2=x_4^2+x_5^2+x_4x_5,~ z_2=x_4x_5(x_4+x_5).
\end{gather*}
Consider the space $\R^4$ with the coordinates $y_1,z_1,y_2,z_2$ and 
the action of the group $\on{S}_2$ on it induced by the action of 
$\on{S}_2=W(U)/(\on{S}_3 \times \on{S}_3)$ on the above basic generators. 
It is easy to check that this action coincides with the diagonal action 
of $\on{S}_2$ on $(\R^2)^2$ for the standard action of $\on{S}_2$ on $\R^2$. 
Since the algebra of $\on{S}_2$-invariant polynomials
on $(\R^2)^2$ is generated by the polarizations of basic invariants for 
the standard action of $\on{S}_2$ ob $\R^2$ we get the following system of 
generators of $\R[U]^{W(U)}$:
\[
f_1=y_1+y_2,~ f_2=z_1+z_2,~ f_3=y_1^2+y_2^2,~ f_4=y_1z_1+y_2z_2,~ f_5=z_1^2+z_2^2.
\]
Simple calculations for the restrictions of the Newton polynomials $N_i$ 
on $\R^6$ to $U$
gives the following result:
\begin{gather*}
N_1|_U=0,~
N_2|_U=2f_1,~
N_3|_U=-3f_2,~\\
N_4|_U=2f_3,~
N_5|_U=-5f_4,~
N_6|_U=3f_5+3f_1f_3-f_1^3.
\end{gather*}
This proves that the morphism $\R[\R^6]^{\on{S}_6}\to\R[U]^{W(U)}$ 
defined by restriction is surjective.
\end{example}

\end{document}